\documentclass[12pt]{article}

\usepackage{amsmath,amssymb,amsthm,amscd}

\begin{document}

\newcommand{\ts}{\,}
\newcommand{\tss}{\hspace{1pt}}
\newcommand{\Mat}{{\rm{Mat}}}
\newcommand{\CC}{\mathbb{C}}
\newcommand{\Sym}{\mathfrak S}

\newtheorem{thm}{Theorem}
\theoremstyle{definition}
\newtheorem{examp}[thm]{Example}
\theoremstyle{remark}
\newtheorem{remark}[thm]{Remark}
\newcommand{\bth}{\begin{thm}}
\renewcommand{\eth}{\end{thm}}
\newcommand{\bex}{\begin{examp}}
\newcommand{\eex}{\end{examp}}
\newcommand{\bre}{\begin{remark}}
\newcommand{\ere}{\end{remark}}

\newcommand{\bal}{\begin{aligned}}
\newcommand{\eal}{\end{aligned}}
\newcommand{\beq}{\begin{equation}}
\newcommand{\eeq}{\end{equation}}
\newcommand{\ben}{\begin{equation*}}
\newcommand{\een}{\end{equation*}}

\newcommand{\bpf}{\begin{proof}}
\newcommand{\epf}{\end{proof}}

\def\beql#1{\begin{equation}\label{#1}}
\title{\Large\bf  On the fusion procedure for the symmetric group}

\author{{\sc A. I. Molev}\\[15pt]
School of Mathematics and Statistics\\
University of Sydney,
NSW 2006, Australia\\
{\tt
alexm\hspace{0.09em}@\hspace{0.1em}maths.usyd.edu.au
}
}

\date{} 

\maketitle

\begin{abstract}

We give a new version of the fusion procedure
for the symmetric group which originated in the work of Jucys
and was developed by Cherednik.
We derive it from the Jucys--Murphy formulas for the diagonal
matrix units for the symmetric group.
\end{abstract}

\bigskip


\section{Introduction}

A key role in the quantum inverse scattering method is played
by solutions of the Yang--Baxter equation.
The {\it fusion procedure\/} is commonly understood as a certain way
to obtain new solutions of this equation
out of the old ones. Consider the Yang--Baxter equation
with the spectral parameter
\beql{ybe}
R_{12}(u)\ts R_{13}(u+v)\ts R_{23}(v)=
R_{23}(v)\ts R_{13}(u+v)\ts R_{12}(u),
\eeq
where $R(u)$ is a function of $u$ with values in the
endomorphism algebra of the tensor square $V\otimes V$
of a vector space $V$.
Both sides of \eqref{ybe}
are endomorphisms of $V\otimes V\otimes V$ and the subscripts
of $R(u)$
indicate the copies of $V$ so that $R_{12}(u)=R(u)\otimes 1$, etc.
The Yang $R$-matrix is
a simplest solution of \eqref{ybe}, which is given by
\beql{rmatrix}
R(u)=1-P\tss u^{-1},
\eeq
where $P$ is the permutation operator
\ben
P:\xi\otimes\eta \mapsto \eta\otimes\xi,\qquad \xi,\eta\in V.
\een
Note that for the values $u=-1$ or $u=1$ the endomorphism
\eqref{rmatrix} maps the space $V\otimes V$ into the subspace
of symmetric or anti-symmetric tensors, respectively.
For instance, take $W$ to be the subspace of symmetric tensors.
Due to \eqref{ybe}, the subspace $V\otimes W$ of $V\otimes V\otimes V$
is preserved by the operator $R_{12}(u)\ts R_{13}(u-1)$.
Similarly, introducing an extra copy of $V$ labeled by $0$,
we can define the operator
\beql{rv}
R_{\tss W}(u)=
R_{12}(u+1)\ts R_{13}(u)\ts
R_{02}(u)
\ts R_{03}(u-1).
\eeq
Then the restriction
of $R_{\tss W}(u)$ to the subspace $W\otimes W$ is
a ``fused $R$-matrix" which is again a solution
of the Yang--Baxter equation.

More generally, for any
different complex numbers $c_1,\dots,c_n$
the operator
\ben
R_{01}(u-c_1)\tss R_{02}(u-c_2)\dots R_{0n}(u-c_n)
\een
preserves the subspace $V\otimes W\subseteq V\otimes V^{\otimes n}$,
where
\beql{wim}
W=\Big(\prod_{1\leqslant i<j\leqslant n}
R_{ij}(c_i-c_j)\Big)\ts V^{\otimes n},
\eeq
and the factors are taken in the lexicographical ordering
on the pairs $(i,j)$.

The symmetric group $\Sym_n$ acts naturally on $V^{\otimes n}$
by permutations of the tensor factors.
When the parameters $c_i$ are chosen in a certain particular way,
the product of the $R$-matrices in \eqref{wim} turns out to
coincide with the image of a diagonal matrix element
of an irreducible representation of $\Sym_n$.
Hence, due to the Schur--Weyl duality,
this leads, by analogy with the above example,
to the construction of solutions of the Yang--Baxter
equation of the type $R_{W}(u)$, where $W$ is an arbitrary
polynomial representation of $GL_N$.
More precisely, the $c_i$ should be taken to be equal
to the respective {\it contents\/}
of a standard tableau $T$ associated with a partition of $n$.
It may happen, however, that $c_i=c_j$ for some $i<j$ so that
the corresponding factor $R_{ij}(c_i-c_j)$ in \eqref{wim}
is not defined. Nevertheless, it turns out that
the product can be interpreted as a well-defined operator
via a certain limiting procedure. Such a procedure
providing an expression
for the matrix elements
of irreducible representations of $\Sym_n$ originates
in the work of Jucys~\cite{j:yo}. A similar approach was
developed by Cherednik~\cite{c:sb} in greater generality
for representations of Hecke algebras.
Cherednik's paper does not contain
complete proofs, however. More details were given by
Jimbo, Kuniba, Miwa and Okado~\cite[Lemmas~3.2 and A.1]{jkmo:fm} while
a complete proof of a version of the fusion theorem
was given by Nazarov~\cite[Theorem~2.2]{n:yc}
with simpler arguments than in \cite{j:yo};
see also Guizzi and Papi~\cite{gp:ca}.
Nazarov's theorem establishes a continuity property
of the restriction of the product of the $R$-matrices
in \eqref{wim} on a certain subset of the parameters $c_i$.
A hook version of the fusion procedure was developed in
recent work of Grime~\cite{g:hf}.

In this paper we give a new version of the fusion procedure
which is similar to \cite{j:yo} but with a different definition
of the limits of rational functions: we take them
{\it consecutively\/} for each single variable.
We show that the procedure is essentially equivalent
to another construction of
Jucys~\cite{j:fy} (which was re-discovered by Murphy~\cite{m:nc})
providing explicit formulas for the
diagonal matrix elements of irreducible representations of $\Sym_n$
in terms of certain elements of the group
algebra $\CC[\Sym_n]$, known as the {\it Jucys--Murphy elements\/}.
The proof of these formulas is rather simple
(its version is reproduced below) thus leading to a short
derivation of the fusion procedure.
Comparing this derivation with the other proofs, note that
the approaches of Jucys~\cite{j:yo} and Nazarov~\cite{n:yc}
rely on the formulas for the diagonal matrix
elements involving the Young symmetrizers and do not
establish a direct relationship with the Jucys--Murphy
construction. Some versions of the procedure
were given by Nazarov~\cite{n:mh} and Grime~\cite{g:hfq}
for the Hecke algebra, and Nazarov~\cite{n:rr} developed
a ``skew fusion procedure";  see also earlier results of
Cherednik~\cite{c:sb, c:ni}.

In what follows, we will work with the group algebra $\CC[\Sym_n]$,
since the product of the $R$-matrices in \eqref{wim}
obviously coincides with the image of the ordered product
\ben
\prod_{1\leqslant i<j\leqslant n}\varphi_{ij}(c_i,c_j),\qquad
\varphi_{ij}(u,v)=1-\frac{(i\ts j)}{u-v}.
\een

\section{Young basis}
\label{sec:yb}

Let us fix some notation and recall some well known facts
about the representations of the symmetric group $\Sym_n$; see e.g.
\cite{jk:rt}.
We write a partition $\lambda$ as a sequence
$\lambda=(\lambda_1,\dots,\lambda_l)$ of integers such that
$\lambda_1\geqslant\dots\geqslant\lambda_l\geqslant 0$.
We shall identify
a partition $\lambda$ with its
diagram which is
a left-justified array of rows of cells such that the top row
contains $\lambda_1$ cells, the next row contains
$\lambda_2$ cells, etc.
Let us fix a positive integer $n$.
If $\lambda_1+\dots+\lambda_l=n$ then $\lambda$ is a partition of $n$,
written $\lambda\vdash n$.
A cell of $\lambda$ is called {\it removable\/} if its removal
leaves a diagram. Similarly, a cell is {\it addable\/} to $\lambda$
if the union of $\lambda$ and the cell is a diagram.
We shall write $\mu\to\lambda$ if $\lambda$ is obtained from $\mu$
by adding one cell.
A tableau $T$ of shape $\lambda$ (or a $\lambda$-tableau $T$)
is obtained by
filling in the cells of the diagram
bijectively with the numbers $\{1,\dots,n\}$.
We write $\text{sh}(T)=\lambda$ if the shape of $T$ is $\lambda$.
A tableau $T$ is called standard if its entries
strictly increase along the rows and down the columns.

The irreducible representations of $\Sym_n$ over $\CC$
are parameterized by partitions of $n$.
Given a partition $\lambda$ of $n$ denote the corresponding irreducible
representation of $\Sym_n$ by $V_{\lambda}$.
The vector space $V_{\lambda}$ is equipped with an $\Sym_n$-invariant
inner product $(\ ,\ )$.
The orthonormal Young basis
$\{v_T\}$ of $V_{\lambda}$ is
parameterized by the set of standard $\lambda$-tableaux $T$.
The action of the standard generators $s_i=(i,i+1)$
of $\Sym_n$ in the Young basis is described as follows.
If $\alpha$ is a cell of $\lambda$ which occurs in row $i$ and column $j$
then the {\it content\/} of $\alpha$ is the number $j-i$.
Now let a standard tableau $T$ be given.
We denote by $c_k=c_k(T)$ the content
of the cell occupied by the number $k$.
Then for any $i\in\{1,\dots,n-1\}$ we have
\ben
s_i\cdot v_T=d\tss v_T+\sqrt{1-d^2}\ts v_{s_iT},
\een
where $d=(c_{i+1}-c_i)^{-1}$, the
tableau $s_iT$ is obtained from $T$ by swapping the entries
$i$ and $i+1$, and
we assume $v_{s_iT}=0$ if the tableau $s_iT$ is not standard.

The group algebra $\CC[\Sym_n]$ is isomorphic to the direct sum
of matrix algebras
\ben
\CC[\Sym_n]\cong \underset{\lambda\vdash n}\bigoplus\ts\ts
\Mat_{f_{\lambda}}(\CC),
\een
where $f_{\lambda}=\dim V_{\lambda}$. The matrix units
$E_{TT'}\in \Mat_{f_{\lambda}}(\CC)$ are parameterized by pairs
of standard $\lambda$-tableaux $T$ and $T'$.
An isomorphism between the algebras is
provided by the formulas
\ben
E_{TT'}=\frac{f_{\lambda}}{n!}\ts\Phi_{TT'},
\een
where $\Phi_{TT'}$ is the matrix element corresponding to the
basis vectors $v_T$ and $v_{T'}$ of the representation $V_{\lambda}$,
\ben
\Phi_{TT'}=\sum_{s\in \Sym_n}(s\cdot v^{}_{T},v^{}_{T'})\cdot s^{-1}
\in\CC[\Sym_n].
\een
In what follows we only use the diagonal matrix units so we
shall write $E_{T}=E_{TT}$ and $\Phi_{T}=\Phi_{TT}$.
Now recall the construction of the matrix units $E_T$
which is due to Jucys~\cite{j:fy} and Murphy~\cite{m:nc}.
Consider the {\it Jucys--Murphy\/} elements $X_1,\dots,X_{n}\in\CC[\Sym_n]$
given by
\ben
X_1=0,\qquad X_i=(1\ts i)+(2\ts i)+\dots+(i-1\ts i),\quad i=2,\dots,n.
\een
The vectors of the Young basis are eigenvectors
for the action of $X_i$ on $V_{\lambda}$:
for any standard $\lambda$-tableau $T$ we have
\beql{xivt}
X_i\cdot v_{T}=c_i(T)\ts v_{T},\qquad i=2,\dots,n.
\eeq
For any $n\geqslant 2$ we regard $\Sym_{n-1}$ as the natural subgroup
of $\Sym_n$.
The branching properties of the Young basis imply the following
properties of the matrix units.
If $U$
is a given standard tableau with the entries $1,\dots,n-1$ then
\beql{eutet}
E_U=\sum_{U\to T} E_T,
\eeq
where $U\to T$ means that
the standard tableau $T$ is obtained from $U$
by adding one cell with the entry $n$.
Relations \eqref{xivt} imply
\beql{xiet}
X_i\ts E_T=E_T\ts X_i=c_i(T)\ts E_T, \qquad i=2,\dots,n
\eeq
for any standard $\lambda$-tableau $T$. In particular,
we have the identity in $\CC[\Sym_n]$,
\beql{xncnet}
X_n=\sum_{\lambda\vdash n}\sum_{\text{sh}(T)=\lambda} c_n(T)\ts E_T.
\eeq
Obviously, $E_{T_0}=1$
if $T_0$ is the $(1)$-tableau with the entry $1$.
The other matrix units are given by the following recurrence relation
which yields an explicit expression of $E_T$ in terms of the
Jucys--Murphy elements $X_2,\dots,X_{n}$.
Let $\lambda\vdash n$ for $n\geqslant 2$ and let $T$ be
a standard $\lambda$-tableau. Let $U$ be the standard tableau
obtained from $T$ by removing the cell $\alpha$ occupied by $n$ and
let $\mu$ be the shape of $U$. Then
\beql{murphyfo}
E_T=E_U\ts \frac{(X_n-a_1)\dots (X_n-a_k)}{(c-a_1)\dots (c-a_k)},
\eeq
where $a_1,\dots,a_k$ are the contents of all addable cells of $\mu$
except for $\alpha$, while $c$ is the content of the latter.
The relation follows from \eqref{eutet} and \eqref{xiet}.
It admits the following interpretation.
Let $u$ be a complex variable. Due to \eqref{xncnet}, the following
is a well-defined rational function
in $u$ with values in $\CC[\Sym_n]$,
\ben
E_T(u)=E_U\ts \frac{u-c}{u-X_n}.
\een
Then $E_T(u)$ is regular
at $u=c$ and $E_T(c)=E_T$.
Indeed,
by \eqref{eutet} and \eqref{xiet} we have
\ben
E^{}_U\ts \frac{u-c\ }{u-X_n}=\sum_{U\to T'}
E^{}_{T'}\ts\frac{u-c\ }{u-c_n(T')}
=E^{}_{T}+\sum_{U\to T',\ts T'\ne T}
E^{}_{T'}\ts\frac{u-c\ }{u-c_n(T')}.
\een
Since $c_n(T')\ne c$ for all standard tableaux $T'$
distinct from $T$,
the value of this rational function at $u=c$ is $E^{}_{T}$.
Thus, the Jucys--Murphy formula \eqref{murphyfo} can also be written as
\ben
E_T=E_U\ts \frac{u-c}{u-X_n}\ts\Big|^{}_{u=c}.
\een
We shall need the corresponding relation for the
matrix elements $\Phi_U$ and $\Phi_T$. Recalling that
the ratio $n!/f_{\lambda}$ equals the product of the hooks of $\lambda$,
we get
\beql{phith}
\Phi_T=
H_{\lambda,\mu}\ts
\Phi_U\ts \frac{u-c}{u-X_n}\ts\Big|^{}_{u=c},
\eeq
where the coefficient $H_{\lambda,\mu}$
is the ratio of the product of hooks of $\lambda$ and
the product of hooks of $\mu$. It can be
given by
\beql{hlamu}
H_{\lambda,\mu}=\frac{(a_1-c)\dots (a_p-c)(c-a_{p+1})\dots (c-a_k)}
{(b_1-c)\dots (b_q-c)(c-b_{q+1})\dots (c-b_r)},
\eeq
where the numbers $a_1,\dots,a_p,c,a_{p+1},\dots,a_k$ are the contents of
all addable cells of $\mu$
and $b_1,\dots,b_q,c,b_{q+1},\dots,b_r$ are the contents
of all removable cells of $\lambda$
with both sequences written in the decreasing order.

\medskip
\noindent
{\it Remark.}
Consider the character $\chi_{\lambda}$ of $V_{\lambda}$,
\ben
\chi_{\lambda}=\sum_{s\in \Sym_n}\chi_{\lambda}(s)\ts s\in \CC[\Sym_n].
\een
We have
\ben
\chi_{\lambda}=\sum_T\Phi_T,
\een
summed over all standard $\lambda$-tableaux $T$. Formula~\eqref{murphyfo}
implies a recurrence relation
for the normalized characters
$\widehat\chi_{\lambda}=f_{\lambda}\ts\chi_{\lambda}/{n!}$,
\ben
\widehat\chi_{\lambda}=\sum_{\mu\to\lambda}\widehat\chi_{\mu}
\frac{(X_n-a_1)\dots (X_n-a_k)}{(c-a_1)\dots (c-a_k)}.
\een
Equivalently,
\ben
\chi_{\lambda}=\sum_{\mu\to\lambda}\chi_{\mu}
\frac{(a_1-X_n)\dots (a_p-X_n)(X_n-a_{p+1})\dots (X_n-a_k)}
{(b_1-c)\dots (b_q-c)(c-b_{q+1})\dots (c-b_r)},
\een
with the notation used in \eqref{hlamu}.

\section{Fusion procedure}
\label{sec:fp}

For any distinct indices $i,j\in\{1,\dots,n\}$ introduce
the rational function in two variables $u,v$
with values in the group algebra $\CC[\Sym_n]$ by
\ben
\varphi_{ij}(u,v)=1-\frac{(i\ts j)}{u-v}.
\een
Equip the set of all pairs $(i,j)$ with $1\leqslant i<j\leqslant n$
with the reverse lexicographical ordering so that $(i_1,j_1)$
precedes $(i_2,j_2)$ if $j_1<j_2$ or $j_1=j_2$ and $i_1<i_2$.
Take
$n$ complex variables $u_1,\dots,u_n$ and consider the ordered product
\ben
\Phi(u_1,\dots,u_n)=
\prod_{1\leqslant i<j\leqslant n}^{\longrightarrow}\varphi_{ij}(u_i,u_j).
\een
Note that the product taken in the (direct) lexicographical ordering
on the pairs $(i,j)$ defines the same rational function.
Now let $\lambda\vdash n$ and fix a standard $\lambda$-tableau $T$.
Set $c_i=c_i(T)$ for $i=1,\dots,n$.

\medskip
\noindent
{\bf Theorem.\ }{\sl
The consecutive evaluations
\ben
\Phi(u_1,\dots,u_n)\big|_{u_1=c_1}\big|_{u_2=c_2}\dots \big|_{u_n=c_n}
\een
of the rational function $\Phi(u_1,\dots,u_n)$ are well-defined.
The corresponding value coincides with
the matrix element $\Phi_T$.
}

\bpf
Clearly, it is sufficient to consider
the last evaluation $u_n=c_n$.
We argue by induction on $n$
and suppose that $n\geqslant 2$.
By the induction hypothesis, setting $u=u_n$ we get
\ben
\Phi(u_1,\dots,u_n)\big|_{u_1=c_1}\dots \big|_{u_{n-1}=c_{n-1}}
=\Phi_U\ts
\varphi_{1n}(c_1,u)\dots \varphi_{n-1,n}(c_{n-1},u),
\een
where the standard tableau $U$ is obtained from $T$ by
removing the cell occupied by $n$. Let us verify that the expression
on the right hand side can be given by
\beql{phiupr}
\Phi_U\ts
\varphi_{1n}(c_1,u)\dots \varphi_{n-1,n}(c_{n-1},u)
=\prod_{i=1}^{n-1}\Big(1-\frac{1}{(u-c_i)^2}\Big)\ts\Phi_U\ts
(1-X_n\ts u^{-1})^{-1}.
\eeq
Note that due to \eqref{xncnet}, the expression
$(1-X_n\ts u^{-1})^{-1}$ is a well-defined rational function in $u$.
Since
\ben
\varphi_{in}(c_i,u)^{-1}\Big(1-\frac{1}{(u-c_i)^2}\Big)=
\varphi_{in}(-c_i,-u),
\een
relation \eqref{phiupr} is equivalent to
\ben
\Phi_U\ts \varphi_{n-1,n}(-c_{n-1},-u)\dots \varphi_{1n}(-c_1,-u)=
\Phi_U\ts
(1-X_n\ts u^{-1}).
\een
We verify by induction on $n$
a slightly more general identity
\begin{multline}\label{phiuinge}
\Phi_U\ts \varphi_{n-1,r}(-c_{n-1},-u)\dots \varphi_{1r}(-c_1,-u)\\
=
\Phi_U\ts
\Big(1-\frac{(1\ts r)+(2\ts r)+\dots+(n-1\ts r)}{u}\ts\Big),
\end{multline}
where $r$ is a fixed index, $r\geqslant n$.
By \eqref{eutet} we
can write $\Phi_U$ as the product
\ben
\Phi_U=\gamma\cdot\Phi_U\ts \Phi_Y,
\een
where $Y$ is the standard tableau obtained from $U$ by
removing the cell occupied by $n-1$ and $\gamma$ is a nonzero constant.
Hence, using the induction hypothesis
we can transform the left hand side of \eqref{phiuinge} as
\ben
\bal
&\gamma\cdot\Phi_U\ts \Phi_Y\ts
\varphi_{n-1,r}(-c_{n-1},-u)\dots \varphi_{1r}(-c_1,-u)\\
{}={}&\gamma\cdot\Phi_U\ts
\varphi_{n-1,r}(-c_{n-1},-u)\ts\Phi_Y\ts\varphi_{n-2,r}(-c_{n-2},-u)
\dots \varphi_{1r}(-c_1,-u)\\
{}={}&\gamma\cdot\Phi_U\ts
\varphi_{n-1,r}(-c_{n-1},-u)\ts\Phi_Y\ts
\Big(1-\frac{(1\ts r)+(2\ts r)+\dots+(n-2\ts r)}{u}\ts\Big).
\eal
\een
This equals
\beql{pcfhiu}
\Phi_U\ts\Big(1-\frac{(n-1\ts r)}{u-c_{n-1}}\Big)
\Big(1-\frac{(1\ts r)+(2\ts r)+\dots+(n-2\ts r)}{u}\ts\Big).
\eeq
Now observe that
\ben
(n-1\ts r)\big((1\ts r)+(2\ts r)+\dots+(n-2\ts r)\big)=
X_{n-1}\ts (n-1\ts r)
\een
and recall that $\Phi_U\ts X_{n-1}=c_{n-1}\ts \Phi_U$ by \eqref{xiet}.
Hence, \eqref{pcfhiu} simplifies to \eqref{phiuinge} as required.

Now write the right hand side of
\eqref{phiupr} as
\beql{rhsfin}
\prod_{i=1}^{n-1}\Big(1-\frac{1}{(u-c_i)^2}\Big)\ts \frac{u}{u-c_n}
\cdot \Phi_U\ts \frac{u-c_n}{u-X_n}.
\eeq
Observe that the product
\ben
\prod_{i=1}^{n-1}\Big(1-\frac{1}{(u-c_i)^2}\Big)\ts \frac{u}{u-c_n}
\een
only depends on the shape $\mu$ of $U$ so we may choose a particular
(e.g. row-standard) tableau $U$ for its evaluation.
A short calculation shows that this product is regular at $u=c_n$
with the value $H_{\lambda,\mu}$ for $c=c_n$.
Due to \eqref{phith},
the value of \eqref{rhsfin} at $u=c_n$ is $\Phi_T$.
\epf

\smallskip
\noindent
{\it Example.}
Let $\lambda=(2^2)$ so that $n=4$. Take the standard $\lambda$-tableau


\setlength{\unitlength}{0.7em}
\begin{center}
\begin{picture}(16,5)

\put(4,1.6){$T=$}

\put(8.7,2.5){1}
\put(8.7,0.5){3}
\put(10.7,2.5){2}
\put(10.7,0.5){4}

\put(8,4){\line(1,0){4}}
\put(8,2){\line(1,0){4}}
\put(8,0){\line(1,0){4}}

\put(8,0){\line(0,1){4}}
\put(10,0){\line(0,1){4}}
\put(12,0){\line(0,1){4}}

\end{picture}
\end{center}
\setlength{\unitlength}{1pt}

\medskip
\noindent
The contents are $c_1=0$, $c_2=1$, $c_3=-1$, $c_4=0$.
Therefore,
\begin{multline}\nonumber
\Phi(0,1,-1,u)=\Big(1+(1\ts 2)\Big)
\ts\Big(1-(1\ts 3)\Big)
\ts\Big(1-\frac{(2\ts 3)}{2}\Big)\\
{}\times \Big(1+\frac{(1\ts 4)}{u}\Big)
\ts\Big(1+\frac{(2\ts 4)}{u-1}\Big)
\ts\Big(1+\frac{(3\ts 4)}{u+1}\Big).
\end{multline}
By the Theorem, this rational function
is well-defined at $u=0$. The corresponding value
is
\begin{multline}\nonumber
\Phi^{}_{\ts T}=\Phi(0,1,-1,0)
=\Big(1+(1\ts 2)\Big)
\ts\Big(1-(1\ts 3)\Big)
\ts\Big(1-\frac{(2\ts 3)}{2}\Big)\\
{}\times\Big(2-(1\ts 4)-(2\ts 4)-(3\ts 4)\Big)\ts
\Big(2+(1\ts 4)+(2\ts 4)+(3\ts 4)\Big).
\end{multline}

\end{document}